\documentclass[12pt]{article}
\usepackage{latexsym,amsmath,amssymb}

\newfont{\bb}{msbm10 at 12pt}
\newfont{\tbb}{msbm10 at 9pt}
\def\l{\hbox{\bb L}}
\def\r{\hbox{\bb R}}

\def\h{\hbox{\bb H}}
\def\n{\hbox{\bb N}}

\def\s{\hbox{\bb S}}
\def\ts{\hbox{\tbb S}}
\def\th{\hbox{\tbb H}}
\def\tn{\hbox{\tbb N}}

\newcommand{\set}[1]{\left\{#1\right\}}
\newcommand{\meta}[2]{\langle #1,#2 \rangle }

\newcommand{\To}{\longrightarrow }
\newcommand{\Diff}{\mathcal{D}(\s ^n) }
\newcommand{\IsoH}{\mathcal{I}(\h ^{n+1})}
\newcommand{\IsoN}{\mathcal{I}(\n ^{n+1}_+)}
\newcommand{\IsoL}{\mathcal{I}(\l ^{n+2})}
\newcommand{\Th}{T_{\left\vert \th ^{n+1}\right.}}
\newcommand{\Tn}{T_{\left\vert \tn ^{n+1}_+\right.}}

\usepackage[latin1]{inputenc}
\usepackage {amsmath}
\usepackage {amsthm}
\usepackage{times}
\usepackage{amscd}
\usepackage{epsf}

\numberwithin{equation} {section}

\begin{document}
\mbox{}\vspace{0.4cm}\mbox{}

\theoremstyle{plain}\newtheorem{lemma}{Lemma}[section]
\theoremstyle{plain}\newtheorem{proposition}{Proposition}[section]
\theoremstyle{plain}\newtheorem{theorem}{Theorem}[section]
\theoremstyle{plain}\newtheorem{example}{Example}[section]
\theoremstyle{plain}\newtheorem{remark}{Remark}[section]
\theoremstyle{plain}\newtheorem{corollary}{Corollary}[section]
\theoremstyle{plain}\newtheorem{definition}{Definition}[section]

\begin{center}
\rule{15cm}{1.5pt} \vspace{.6cm}

{\Large \bf Invariant conformal metrics on $\s ^n$} \vspace{0.4cm}

\vspace{0.5cm}

{\large José M. Espinar$\,^\dag$\footnote{The author is partially supported by
Spanish MEC-FEDER Grant MTM2007-65249, and Regional J. Andalucia Grant
P06-FQM-01642}}\\
\vspace{0.3cm} \rule{15cm}{1.5pt}
\end{center}

\vspace{.5cm}

\noindent $\mbox{}^\dag$ Departamento de Geometría y Topología, Universidad de
Granada, 18071 Granada, Spain; e-mail: jespinar@ugr.es\vspace{0.2cm}

\vspace{.3cm}


\vspace{.3cm}

\begin{abstract}
In this paper we use the relationship between conformal metrics on the sphere and
horospherically convex hypersurfaces in the hyperbolic space for giving sufficient
conditions on a conformal metric to be radial under some constrain on the
eigenvalues of its Schouten tensor. Also, we study conformal metrics on the sphere
which are invariant by a $k-$parameter subgroup of conformal diffeomorphisms of the
sphere, giving a bound on its maximum dimension.

Moreover, we classify conformal metrics on the sphere whose eigenvalues of the
Shouten tensor are all constant (we call them \emph{isoparametric conformal
metrics}), and we use a classification result for radial conformal metrics which are
solution of some $\sigma _k -$Yamabe type problem for obtaining existence of
rotational spheres and Delaunay-type hypersurfaces for some classes of Weingarten
hypersurfaces in $\h ^{n+1}$.
\end{abstract}

\section{Introduction}

In the last 30 years, the \emph{Nirenberg Problem}, i.e., \emph{which functions
$S:\s ^n \To \r$ arise as the scalar curvature of some conformal metric on the
sphere?}, has received an amazing number of contributions (see
\cite{AM,BC,BE,Ch,CY1,CY2,CLi1,CLi2,CLn,ES,KW1,KW2,Li1,Li2,Mo}), but sufficient and
necessary conditions for the solvability are still unknown.

However, this problem opened the door of a rich subject in the last few years,
\emph{conformally invariant equations}. Let $\mathcal{F}(x_1 , \ldots , x_n)$ denote
a smooth functional, and let $\Gamma \in C^{\infty}(\s ^n)$. \emph{Does there exist
a conformal metric $g = e^{2\rho}g_0$ on $\s ^n$ such that the eigenvalues $\lambda
_i$ of its Schouten tensor verify
$$ \mathcal{F}(\lambda _1 , \ldots , \lambda _n ) = \Gamma , \, \text{ on } \s ^n .$$}

Given $(M, g)$ a Riemannian manifold, for $n \geq 3$, the Schouten tensor of $g$ is
given by
$$ {\rm Sch}_g := \frac{1}{n-2}\left( {\rm Ric}(g) -
\frac{S(g)}{2(n-1)}g\right)$$where ${\rm Ric}(g)$ and $S(g)$ are the Ricci tensor
and the scalar curvature function of $g$ respectively.

Note that, when $\mathcal{F}(x_1 , \ldots , x_n ) = x_1 +\cdots + x_n$ we have the
Nirenberg Problem. Right now, the most developed topic for these equations is when
we consider $\mathcal{F}(\lambda _1 , \ldots , \lambda _n)\equiv \sigma _k (\lambda
_i)$ as the $k-$th elementary symmetric polynomial of its arguments equals to a
constant, i.e.,
\begin{equation}\label{1}
\sigma _k (\lambda _i) = {\rm constant}.
\end{equation}

Many deep results are known for these equations (see
\cite{Ch,CGY1,CGY2,Li3,Li4,LL1,LL2,Vi} and reference therein). Mostly of these
results are devoted to solutions either on $\s ^n$ or $\r ^n$, and little is known
when we look for conformal metrics on a domain of the sphere (see \cite{LN,MP} and
references therein). In this line, Chang-Han-Yang \cite{CHY} have classified all
posible radial solution to the equation \eqref{1} \emph{``as guidance in studying
the behavior of singular solutions in the general situation''}. This is natural
since radial solutions are the simplest examples. Thus, the next step is:
\emph{under what (local) conditions can we know that the solution is radial?}

In a recent paper \cite{EGM}, the authors showed a correspondence between conformal
metrics on the sphere and horospherically convex hypersurfaces in the hyperbolic
space. Here, they provide a back-and-forth construction which give a hypersurface
theory interpretation for the famous Nirenberg Problem, relating it with a natural
formulation of the Christoffel problem in $\h ^{n+1}$. Moreover, this correspondence
is more general and it relates conformally invariant equations with Weingarten
hypersurfaces horospherically convex. The main line in this paper is to use the deep
theorems on conformal geometry to infer results in the hypersurface theory, but,
\emph{how can the hypersurface theory help to get information on conformal
geometry?}

We will see here that, using the hypersurface setting, we can obtain sufficient
conditions under which a conformal metric is radial among others on invariant
conformal metrics under a subgroup of conformal diffeomorphisms of the sphere. We
should mention that the theorems included here are local results, besides the usuals
on this direction that they are from a global character.

In Section $2$ we establish the necessary preliminaries on conformal geometry, and
it is devoted also to summarize the correspondence developed in \cite{EGM} between
conformal metrics and horospherically convex hypersurfaces, that is, given a
conformal metric on the sphere they construct a horospherically convex hypersurface
in $\h^{n+1}$ and viceversa. In Section $3$ we establish that if a conformal metric
is invariant under a subgroup of conformal diffeomorphism of the sphere, then its
associated horospherically convex hypersurface is invariant under the subgroup of
isometries induce by the subgroup of conformal diffeomorphism, and viceversa, i.e.,

\vspace{.3cm}

{\bf Lemma \ref{confisoinv}:}

{\it Let $\phi: \Omega \subset \s ^n\To \h^{n+1}$ be a locally horospherically
convex hypersurface with hyperbolic Gauss map $G(x)=x$, support function
$e^{\rho}:\Omega \To (0,+\infty)$, and let $g=e^{2\rho}g_0$ denote its horospherical
metric. Let $\Th \in \IsoH $ be an isometry and $\Phi \in \Diff$ its associated
conformal diffeomorphism. Thus, if $\phi $ is $\Th -$invariant then $g$ is $\Phi
-$invariant.

Conversely, let $g= e^{2\rho} g_0$ be a conformal metric defined on a domain of the
sphere $\Omega \subset \s ^n$ such that the eigenvalues of its Schouten tensor,
${\rm Sch}_g$, verify
$${\rm sup}\set{\lambda _i (x) , \, x \in \Omega , \, i = 1 ,\ldots , n }< +\infty .$$

Let $\Phi \in \Diff $ be a conformal difeomorphism and $\Th \in \IsoH$ its
associated isometry. Thus, if $g $ is $\Phi-$invariant then $\phi$, given by
\eqref{repfor}, is $\Th -$invariant.}

\vspace{.3cm}

In Section $4$ we classify the conformal metrics on the sphere whose eigenvalues of
its Shouten tensor are all constant, we call these metrics \emph{isoparametric
conformal metrics}. Since the above classification have not been done before (as far
as we know), we will include here.

In Section $5$ we state our main results, we give sufficient conditions under which
a conformal metric is radial in terms of the eigenvalues of its Shouten tensor,

\vspace{.3cm}

{\bf Theorem \ref{t1}:}

{\it Let $g= e^{2\rho} g_0$ be a conformal metric defined on a domain of the sphere
$\Omega \subset \s ^n$ such that the eigenvalues, $\lambda_i$, for $i =1, \ldots,
n$, of its Schouten tensor, ${\rm Sch}_g$, verify
$${\rm sup}\set{\lambda _i (x) , \, x \in \Omega , \, i = 1 ,\ldots , n }< +\infty .$$

Furthermore, assume that the eigenvalues satisfy
\begin{equation*}
\begin{split}
\lambda &= \lambda _1 = \cdots = \lambda _{n-1} \\
\nu = \nu (\lambda) &= \lambda _n \\
\lambda - \nu & \neq 0 .
\end{split}
\end{equation*}

Then, $g$ is radial.}

\vspace{.3cm}

Moreover, we study conformal metrics on the sphere which are invariant by a
$k-$parameter subgroup of conformal diffeomorphisms of the sphere, giving a bound on
its maximum dimension,

\vspace{.3cm}

{\bf Theorem \ref{t2}:}

{\it Let $g= e^{2\rho} g_0$ be a conformal metric defined on a domain of the sphere
$\Omega \subset \s ^n$ such that $g \not\in \mathcal{C}(n)$ and the eigenvalues,
$\lambda_i$, for $i =1, \ldots, n$, of its Schouten tensor, ${\rm Sch}_g$, verify
$${\rm sup}\set{\lambda _i (x) , \, x \in \Omega , \, i = 1 ,\ldots , n }< +\infty .$$

Suppose that $g$ is invariant by a $k-$parameter subgroup of conformal
diffeormorphism $\mathcal{G} \leq \Diff$. Then the maximum value of $k$ is $k_{max}=
\frac{n(n-1)}{2}$, and if $k = k_{max}$, the Schouten tensor of $g$, ${\rm Sch}_g$
has two eigenvalues $\lambda$ and $\nu$, where one of them, say $\lambda$, has
multiplicity at least $n-1$. If, in addition, $\lambda \neq 0$, $\nu = \nu
(\lambda)$ and $\nu - \lambda \neq 0$, then $g$ is radial.}

\vspace{.3cm}

Finally, in Section $6$, we give some existence results for some classes of
Weingarten hypersurfaces which are rotationally invariants and horospherically
convex, based on a result of Chang-Han-Yang \cite{CHY}.

\section{Preliminaries}

\subsection{On conformal geometry}

Let $(M ^n , g)$, $n\geq 3$, be a Riemannian manifold. The Riemann curvature tensor,
${\rm Riem}$, can be decomposed as
$$ {\rm Riem} = W_g + {\rm Sch}_g \odot g , $$where $W_g$ is the Weyl tensor,
$\odot$ is the Kulkarni-Nomizu product, and
$$ {\rm Sch}_g := \frac{1}{n-2}\left( {\rm Ric}_g - \frac{S(g)}{2(n-1)}g\right) $$is
the \emph{Schouten tensor}. Here ${\rm Ric}_g$ and $S(g)$ stand for the Ricci
curvature and scalar curvature of $g$ respectively.

The eigenvalues of ${\rm Sch}_g$ are defined as the eigenvalues of the endomorphism
$g^{-1}{\rm Sch}_g$ and we will denote them by $\lambda _i$, $i =1, \ldots , n$.

It is well known that the Schouten tensor encodes all the information on how
curvature varies by a conformal change of metric. It is worth it to remark that the
Weyl curvature tensor vanishes identically when $(M^n ,g)$ is locally conformally
flat since it is the situation of the present work. We will consider conformal
metrics to the standard metric on the $n-$sphere, $(\s ^n ,g_0)$, i.e.,
$$ g = e^{2 \rho} g_0 .$$

\begin{definition}\label{confinv}
Let us denote by $\Diff$ the group of conformal diffeomorphisms on the sphere and
$\Phi \in \Diff$ a conformal diffeomorphism. Let $g = e^{2 \rho} g_0$ be a conformal
metric defined on a domain $\Omega \subset \s ^n$. $g$ is $\Phi -$invariant if
$$ g_x (u,v) = (\Phi ^* g)_x(u,v) , \, \forall x \in \Omega ,  \forall u,v \in T_x \s ^n ,
\text{ such that } \Phi (x) \in \Omega .$$

Moreover, given a continuous subgroup of conformal diffeomorphisms $\mathcal{G} \leq
\Diff$, $g$ is $\mathcal{G}-$invariant if it is $\Phi -$invariant for all $\Phi \in
\mathcal{G}$.
\end{definition}

The basic example of $\mathcal{G}-$invariant metric is that which is radial
symmetric, i.e., when $\mathcal{G}$ is a subgroup of rotations. In this case, we say
that $g$ is radial.

\subsection{On hypersurface theory}

First, let us establish the necessary notation that we will use along the work.
Actually, we will summarize here the construction developed in \cite{EGM} for the
sake of completeness, that is, in order to prove our results, we will use the
correspondence between conformal metrics on the sphere and locally horospherically
convex hypersurfaces in $\h ^{n+1}$. So, we will remind, in a short way, how to
construct a locally horospherically convex hypersurface from a conformal metric on
the sphere.

Let us denote by $\l ^{n+2}$ the $(n+2)-$dimensional Lorentz-Minkowski space, i.e.,
the vectorial space $\r ^{n+2}$ endowed with the Lorentzian metric $\meta{}{}$ given
by

$$ \meta{\bar{x}}{\bar{x}} = - x_0 ^2 + \sum _{i=1}^{n+1} x_i ^2 ,$$where
$\bar{x} \equiv (x_0 , x_1 , \ldots , x_{n+1})\in \r ^{n+2}$.

So, the $(n+1)-$dimensional hyperbolic space, de-Sitter space and null cone are
given, respectively, by the hyperquadrics
\begin{equation*}
\begin{split}
\h ^{n+1} &= \set{ \bar{x} \in \l ^{n+2} : \, \meta{\bar{x}}{\bar{x}} = -1, \, x_0 >0}\\
\s ^{n+1}_1 &= \set{ \bar{x} \in \l ^{n+2} : \, \meta{\bar{x}}{\bar{x}} = 1}\\
\n ^{n+1}_+ &= \set{ \bar{x} \in \l ^{n+2} : \, \meta{\bar{x}}{\bar{x}} = 0, \, x_0
>0} .
\end{split}
\end{equation*}

It is well know that $\h ^{n+1}$ inherits from $(\l ^{n+2}, \meta{}{})$ a Riemannian
metric which make it the standard model of Riemannian space of constant sectional
curvature $-1$. Its ideal boundary at infinity, $\partial _{\infty} \h ^{n+1}$, will
be denoted by $\s ^n _{\infty}$.

Horospheres will play an essential role in what follows, so, we go through
describing their most important properties. In this model, horospheres in $\h
^{n+1}$ are the intersection of affine degenerate hyperplanes of $\l ^{n+2}$ with
$\h ^{n+1}$. Thus, it is clear that the boundary at infinity is a single point $x\in
\s ^n _{\infty}$. In this way, two horospheres are always congruent, and they are at
a constant (hyperbolic) distance if their respective points at infinity agree.
Moreover, given a point $x \in  \s ^{n}_{\infty}$, horospheres having $x$ as its
point at infinity provide a foliation of $\h ^{n+1}$.

From now on, $\phi : M ^n \To \h ^{n+1}$  will denote an oriented immersed
hypersurface and $\eta : M^n \To \s ^{n+1}_1$ its unit normal.

\begin{definition}[\cite{Eps1,Eps2,Br}]
Let $\phi:M^n\To \h ^{n+1}$ denote an immersed oriented hypersurface in $\h^{n+1}$
with unit normal $\eta$. The \emph{hyperbolic Gauss map} $$G:M^n\To \s_{\infty}^n
\equiv \s ^n$$ of $\phi$ is defined as follows: for every $p\in M^n$, $G(p)\in \s
_{\infty}^n$ is the point at infinity of the unique horosphere in $\h ^{n+1}$
passing through $\phi(p)$ and whose inner unit normal at $p$ agrees with $\eta(p)$.
\end{definition}

Associated to $\phi$, let us consider the map
$$ \psi := \phi + \eta : M ^n \To \n ^{n+1}_+ ,$$called the \emph{ associated light cone
map}.

The map $\psi$ is strongly related to the hyperbolic Gauss map $G:M^n\To \s
_{\infty}^n$ of $\phi$. Indeed, the ideal boundary of $\n _+^{n+1}$ coincides with
$\s _{\infty}^n$, and can be identified with the projective quotient space $\n
_+^{n+1} / \r _+$. So, with all of this, we have $G=[\psi] :M^n\To \s
_{\infty}^n\equiv \n _+^{n+1} / \r _+$.

Note that $\psi _0 > 0 $, being $\psi = (\psi _0 , \psi _1 , \ldots , \psi _{n+1})
\in \l ^{n+2}$. If we label $\psi _0 := e^{\rho}$, then we can interpret the
\emph{hyperbolic Gauss map} as the map
$$ G : M ^n \To \s ^n = \n ^{n+1}_{+} \cap \set{x \in \l^{n+2} : \, x_0 =1} $$given
by
\begin{equation}\label{psiG}
\psi = e^{ \rho} (1, G) .
\end{equation}

Moreover, we call $e^{\rho}$ the \emph{horospherical support function}. Also, if
$\{e_1,\dots, e_n\}$ denotes an orthonormal basis of principal directions of $\phi$
at $p$, and if $\kappa _1, \dots, \kappa _n$ are the associated principal
curvatures, it is immediate that

\begin{equation}\label{metpsi}
\meta{(d\psi )_p(e_i)}{(d\psi )_p(e_j)}= (1- \kappa _i)^2 \delta_{ij} = e^{2\rho}
\meta{(dG)_p(e_i)}{(dG )_p(e_j)}_{\ts ^n}.
\end{equation}

Coming back to horospheres, we must remark that horospheres are the unique
hypersurfaces such that, innerly oriented (i.e., when the unit normal points to the
convex side), its associated light cone map is constant: $\phi + \eta =v \in \n
^{n+1}_+$. Moreover, if we write $v = e^{\rho}(1,x)$, we see that $x \in \s ^n$ is
the point at infinity of the horosphere and $\rho$ is the signed hyperbolic distant
of the horosphere to the point $\mathcal{O}=(1,0,\ldots ,0) \in \h ^{n+1}\subset \l
^{n+2}$.

In the hyperbolic setting we have a notion of convexity weaker than the usual
geodesic convexity, i.e.,

\begin{definition}[\cite{Sc}]\label{horocon}
Let $M^n\subset \h^{n+1}$ be an immersed oriented hypersurface, and let
$\mathcal{H}_p$ denote the horosphere in $\h^{n+1}$ that is tangent to $M^n$ at $p$,
and whose interior unit normal at $p$ agrees with the one of $M^n$. We will say that
$M^n$ is \emph{horospherically convex} at $p$ if there exists a neighborhood
$V\subset M^n$ of $p$ so that $V\setminus \{p\}$ does not intersect $\mathcal{H}_p$,
and in addition the distance function of the hypersurface to the horosphere does not
vanish up to the second order at $p$ in any direction.
\end{definition}

Thus, from \eqref{metpsi}, we have the following characterization of horospherically
convex hypersurfaces:

\begin{lemma}\label{hc}
Let $\phi:M^n\To \h^{n+1}$ be an oriented hypersurface. The following conditions are
equivalent at $p\in M^n$.
\begin{enumerate}
 \item[(i)]
All principal curvatures of $M^n$ at $p$ are simultaneously $<1$ or $>1$.
 \item[(ii)]
$M^n$ is horospherically convex at $p$.
 \end{enumerate}

In particular, if $\phi : M^n \To \h^{n+1}$ is horospherically convex at $p$, then
its Gauss map verify $dG_p \neq 0$.
\end{lemma}

So, if $M^n $ is horospherically convex at $p$, therefore $dG_p \neq 0$ and there
exist neighborhoods $U \subset M^n $ and $\Omega \subset \s ^n $ such that $G : U
\To \Omega $ is a diffeomorphism, and
$$ g = e^{2\rho} \meta{dG}{dG}_{\ts ^n} $$define a conformally flat Riemannian metric on $M^n$,
called the \emph{horospherical metric}. Since $G$ is a diffeomorphism between $U$
and $\Omega $ we can use it as a parametrization of the hypersurface, i.e., we can
assume that $\phi : \Omega \subset \s ^n \To \h ^{n+1} $ and $G(x) = x$ on $\Omega
\subset \s ^n$.

Thus, if $G : M ^n \To \Omega \subseteq \s ^n $ is a global diffeormorphism of the
hypersurface onto a domain of the sphere, we can use the hyperbolic Gauss map as a
global parametrization of $\phi$ as above, i.e., $\phi : \Omega \To \h ^{n+1}$ and
$G(x)=x$. In this case, the horospherical metric is given by
$$ g = e^{2\rho} g_0 . $$

Now, we are ready to establish the mentioned relationship between conformal metrics
on the sphere and horospherically convex hypersurfaces

\begin{theorem}[\cite{EGM}]\label{representacion}
Let $\phi: \Omega \subset \s ^n\To \h^{n+1}$ be a horospherically convex
hypersurface with hyperbolic Gauss map $G(x)=x$, support function $e^{\rho}:\Omega
\To (0,+\infty)$, and let $g=e^{2\rho}g_0$ denote its horospherical metric. Then it
holds
\begin{equation}\label{repfor}
\phi = \frac{e^{\rho}}{2}\left( 1+ e^{-2\rho} \left( 1+ ||\nabla^{g_0} \rho
||_{g_0}^2 \right)\right) (1,x) + e^{-\rho} (0, -x +\nabla^{g_0} \rho).
\end{equation}

Moreover, the eigenvalues, $\lambda _i$, of the Schouten tensor of $g$, ${\rm
Sch}_g$, and the principal curvatures, $\kappa _i$, of $\phi$ are related by
\begin{equation}\label{lambdakappa}
\lambda _i = \frac{1}{2} -\frac{1}{1-\kappa _i} .
\end{equation}

Conversely, given a conformal metric $g= e^{2\rho} g_0$ defined on a domain of the
sphere $\Omega \subset \s ^n$ such that the eigenvalues of its Schouten tensor,
${\rm Sch}_g$, are less than $1/2$, then the map $\phi : \Omega \To \h ^{n+1}$ given
by \eqref{repfor} defines a horospherically convex hypersurface in $\h ^{n+1}$ whose
hyperbolic Gauss map is given by $G(x)=x$, $x\in \Omega$.
\end{theorem}

\begin{remark}\label{remark}
We must say that the condition on the eigenvalues of the Schouten tensor is easily
removable, i.e., we only need to ask that
$$ {\rm sup}\set{\lambda _i (x) , \, i= 1, \ldots, n, \, x \in \Omega} < + \infty .$$

If this occurs, we can dilate the metric $g$ as $g_t = e^t g$ for $t>0$. Then, the
eigenvalues of ${\rm Sch}_{g_t}$ are given by
$$ \lambda _i ^{t} = e^{-t} \lambda _i .$$

Thus, for $t$ big enough, we can achieve $\lambda _{i}^t < 1/2$ for $i = 1 , \ldots,
n$.
\end{remark}

\section{Conformal diffeomorphisms and isometries}

Let us denote by $\IsoL $, $\IsoH$ and $\IsoN$ the group of isometries of $\l
^{n+2}$, the $(n+1)-$dimensional hyperbolic space and the $(n+1)-$dimensional null
cone respectively.

It is well known (see \cite{D}) that a conformal diffeormorphism $\Phi \in \Diff$
induce a unique isometry in $\l ^{n+2}$, $T\in \IsoL$, such that restricted to $\h
^{n+1}$ and $\n ^{n+1}_+ $ induces an isometry in these spaces and viceversa. The
restrictions of $T\in \IsoL$ to $\h ^{n+1} $ and $\n ^{n+1}_+$ will be denote by
$\Th$ and $\Tn$ respectively. Moreover, each isometry $\Th \in \IsoH$ induce an
unique isometry $\Tn \in \IsoN $ and viceversa.

\begin{definition}\label{isoinv}
Let $M^n \subset N^n $ be a domain of a $n-$manifold $N$. Let $\phi : M^n \subset
N^n \To \h ^{n+1}$ be a hypersurface and $\Th \in \IsoH$ an isometry. $\phi $ is
$\Th -$invariant if there exists $i_{\Th} : N^n \To N^n $ a diffeomorphism such that
$$ (\Th \circ \phi) (p) = \left(\phi \circ i_{\Th}\right)(p) , \forall p \in M ^n
\text{ such that } i_{\Th}(p)\in M^n .$$

Moreover, given a continuous subgroup of isometries $\mathcal{T} \leq \IsoH$, $\phi$
is $\mathcal{T}-$invariant if it is $\Th -$invariant for all $\Th \in \mathcal{T}$.
\end{definition}

The next result state the relationship between conformal metrics on the sphere which
are invariant by a conformal diffeomorphism and horospherically convex hypersurfaces
which are invariant by an isometry.

\begin{lemma}\label{confisoinv}
Let $\phi: \Omega \subset \s ^n\To \h^{n+1}$ be a locally horospherically convex
hypersurface with hyperbolic Gauss map $G(x)=x$, support function $e^{\rho}:\Omega
\To (0,+\infty)$, and let $g=e^{2\rho}g_0$ denote its horospherical metric. Let $\Th
\in \IsoH $ be an isometry and $\Phi \in \Diff$ its associated conformal
diffeomorphism. Thus, if $\phi $ is $\Th -$invariant then $g$ is $\Phi -$invariant.

Conversely, let $g= e^{2\rho} g_0$ be a conformal metric defined on a domain of the
sphere $\Omega \subset \s ^n$ such that the eigenvalues of its Schouten tensor,
${\rm Sch}_g$, are less than $1/2$. Let $\Phi \in \Diff $ be a conformal
difeomorphism and $\Th \in \IsoH$ its associated isometry. Thus, if $g $ is
$\Phi-$invariant then $\phi$, given by \eqref{repfor}, is $\Th -$invariant.
\end{lemma}
\begin{proof}
On one hand, if $\phi$ is horospherically convex, $\phi $ is $\Th -$invariant if and
only if its associated light cone map $\psi$ is $\Tn -$invariant, i.e., if
\begin{equation}\label{proof}
\left(\Tn \circ \psi \right)(x) = \left( \psi \circ \Phi \right) (x) , \, x \in
\Omega \text{ such that } \Phi (x) \in \Omega .
\end{equation}being $\Phi \in \Diff$ the conformal diffeomorphism associated to $\Tn \in
\IsoN$.

On the other hand, we have an explicit correspondence between conformal
diffeomorphisms on the sphere and isometries on $\n ^{n+1}_+$ (see \cite[Proposition
7.4]{D}). Given an isometry $\Tn \in \IsoN $, at points $(1,x) \in \s ^n = \n
^{n+1}_{+} \cap \set{x \in \l^{n+2} : \, x_0 =1}$ we can see it as
$$ \Tn ((1,x)) = e^{-\omega (x) } (1 , \Phi(x)) ,$$then $ \Phi : \s ^n \To \s ^n$
defines a conformal diffeomorphism on the $n-$sphere with conformal factor
$e^{\omega}$. Conversely, given a conformal diffeomorphism $\Phi \in \Diff $ with
conformal factor $e^{\omega}$, at any point $e^{t}(1,x)\in \n ^{n+1}_+$ define
$$ \Tn (e^{t}(1,x)) = e^t e^{-\omega (x)}( 1, \Phi (x)) ,$$then $\Tn \in \IsoN$.

We first prove the converse. By the previous considerations, we only need to prove
\eqref{proof}. Thus, if $g = e^{2\rho}g_0$ is $\Phi -$invariant, hence by Definition
\ref{confinv} we have that
$$ \rho (x) = \rho (\Phi (x)) + \omega (x), \, \text{ provided } \Phi (x)\in \Omega .$$

Let $\Tn$ be the isometry of $\n ^{n+1}_+$ associated to $\Phi$, then
\begin{equation*}
\begin{split}
\left(\Tn \circ \psi \right) (x) &= \Tn (e^{\rho (x)}(1,x)) = e^{\rho (x) - \omega
(x)}(1,
\Phi (x))\\
 &= e^{\rho (\Phi (x))} (1, \Phi (x)) = \psi (\Phi (x)) = \left( \psi \circ \Phi
 \right)(x).
\end{split}
\end{equation*}

Now, if $\phi $ is $\Th -$invariant, following the above computations, we can
observe that
$$ \rho (x) = \rho (\Phi (x)) + \omega (x), $$being $e^{\omega} :\Omega \to \r$ the
conformal factor of the conformal diffeomorphism, $\Phi$, associated to $\Th$. Thus,
$g$ is $\Phi -$invariant.
\end{proof}

\section{Isoparametric conformal metrics}

Here, we will classify the class of conformal metrics on the sphere such that all
the eigenvalues of its Schouten tensor are constant, we denote this class by
$\mathcal{C}(n)$.

The local classification of conformal metrics on the class $g \in \mathcal{C}(n)$
can be done through a result of E. Cartan \cite{Ca}. Suppose $g \in \mathcal{C}(n)$
therefore, after possibly a dilation, the associated hypersurface given by Theorem
\ref{representacion} is an isoparametric hypersurface in $\h ^{n+1}$, i.e., all its
principal curvatures are constant. Thus, it is a piece of either a totally umbilical
hypersurface (hypersphere, horosphere, totally geodesic hyperplane and equidistant)
or a standard product $\s ^k \times \h ^{n-k}$ in $\h ^{n+1}$. For this reason, we
will call a metric in $\mathcal{C}(n)$ an \emph{isoparametric conformal metric}.

It is known that solutions of
\begin{itemize}
\item $$ \sigma _k (\lambda _i ) = 1  \, \, on \, \, \s ^n  $$are given by conformal
diffeomorphisms of the standard metric on the sphere. Such solution corresponds to a
hypersphere via Theorem \ref{representacion} (see \cite{EGM}).

\item $$ \sigma _k (\lambda _i ) = 0  \, \, on  \, \, \r ^n $$are explictly known (see \cite{Li3}).
Such solution corresponds to a horosphere via Theorem \ref{representacion} (see
\cite{EGM}).
\end{itemize}

Now, our task is to compute explicitly the horospeherical support function
associated to a totally geodesic hyperplane, an equidistant hypersurface and a
standard product $\s ^k \times \h ^{n-k}$. To do so, we will give the
parametrization of such hypersurface and its unit normal vector field and, by means
of equation \eqref{psiG}, we will have an explicit formula for the horospherical
support function and hyperbolic Gauss map. Thus, for an isoparametric hypersurface
$\phi : \Omega \subset \r ^n \To \h ^{n+1} \subset \l ^{n+2}$ with unit normal $\eta
: \Omega \subset \r ^n \To \s _ 1 ^{n+1} \subset \l ^{n+2}$, we will have
$$ \rho : \Omega \subset \r ^n \To \r  ,$$and
$$ G : \Omega \subset \r ^n \To D \subset \s ^{n} \, \, \text{global diffeomorphism}.$$

Hence, the isoparametric conformal metric associated to that hypersurface is given
by
\begin{equation}\label{isop}
g = e^{\rho (G^{-1} (y))} g_0 , \, \, y \in D .
\end{equation}

Let us describe the announced examples:

\begin{enumerate}
\item \textbf{Totally geodesic hyperplanes:}

Set $\Omega = \set{x \in \r ^n : \, |x| < r}$ and $D= \set{y \in \s ^n : \, {\rm
dist}_{g_0}({\bf n},y) < \pi/2 - \arcsin \left(\frac{1-r^2}{1+r^2}\right)}$, where
${\bf n}$ is the north pole. Then,

\begin{eqnarray*}
\psi (x) &=& \left( \frac{1+r^2}{2\sqrt{r^2 -|x|}}, \frac{x}{\sqrt{r^2 - |x|^2}},
\frac{1-r^2}{2\sqrt{r^2 - |x|^2}}\right)\\
\eta (x) &=& \left( \frac{1-r^2}{2r}, {\bf 0}, \frac{1+r^2}{2r}\right) .
\end{eqnarray*}

Thus, from \eqref{psiG}, we get
\begin{eqnarray}
\rho (x) &=& \ln \left( \frac{r(1+r^2)+(1-r^2)\sqrt{r^2 -|x|^2}}{2r\sqrt{r^2 -|x|^2}}\right)\\[3mm]
G (x) \hspace{-3mm} &=&\hspace{-3mm} \left( \frac{2r x}{r(1+r^2)+(1-r^2)\sqrt{r^2
-|x|^2}} , \frac{r(1- r^2)+(1+r^2)\sqrt{r^2 -|x|^2}}{r(1+r^2)+(1-r^2)\sqrt{r^2
-|x|^2}}\right)
\end{eqnarray}

In this case, the principal curvatures all are equals to zero, $k_i =0 $, $i =1,
\ldots, n$. Thus, the eigenvalues of the Shouten tensor associated to $g$ (given by
\eqref{isop}) are $\lambda _i = -1/2$, $i=1,\ldots ,n$.

\item \textbf{Equidistant hypersurfaces:}

Set $\Omega = \set{x \in \r^n : \, |x|<r}$ and $D= \set{y \in \s ^n : \, {\rm
dist}_{g_0}({\bf s},y) < \pi/2 - \arcsin \left(\frac{1-r^2}{1+r^2}\right)}$, where
${\bf s}$ is the south pole. Set $t>0$, $R^2 = t^2 + r^2$ and $\beta (s) = -t +
\sqrt{R^2 -s}$ for $s< r$. Then

\begin{eqnarray*}
\psi (x) &=& \left(\frac{1+|x|^2 +\beta(|x|^2)^2}{2 \beta(|x|^2)},
\frac{x}{\beta(|x|^2)},\frac{1-|x|^2 -\beta (|x|^2)^2}{2 \beta(|x|^2)} \right)\\[3mm]
\eta (x) \hspace{-3mm} &=&\hspace{-3mm}
\left(\frac{(1-t^2-R^2)\sqrt{R^2-|x|^2}+2tR^2}{2R\beta (|x|^2)}, \frac{t x}{R
\beta(|x|^2)},\frac{(1+t^2+R^2)\sqrt{R^2-|x|^2}-2tR^2}{2R\beta (|x|^2)}\right)
\end{eqnarray*}

Thus, from \eqref{psiG}, we get
\begin{eqnarray}
\rho (x) &=& \ln\left(\frac{\alpha (|x|^2)}{2R
\beta(|x|^2)}\right)\\[3mm]
`G (x)\hspace{-2mm} &=&\hspace{-2mm} \left( \frac{2(R-t) x}{\alpha (|x|^2)} ,
\frac{R+ \sqrt{R^2 -|x|^2} + (R+t)^2 (R- \sqrt{R^2 -|x|^2})}{\alpha (|x|^2)}\right)
,
\end{eqnarray}where
$$ \alpha (|x|^2) =R+ \sqrt{R^2 -|x|^2} + (R+t)^2 (R- \sqrt{R^2 -|x|^2}) .$$

In this case, the principal curvatures all are equals to $-t/R$, $k_i =-t/R $, $i
=1, \ldots, n$. Thus, the eigenvalues of the Shouten tensor associated to $g$ (given
by \eqref{isop}) are $\lambda _i = -(R+t)/2(R-t)$, $i=1,\ldots ,n$.

\item {\bf $\h ^{k} \left(-\frac{1}{1+r^2}\right)\times \s ^{n-k} \left( \frac{1}{r}\right)$:}

For the sake of simplicity, we parametrize just a half of this hypersurface.

Set $\Omega = \set{x \in \r ^k } \times \set{z \in \r ^{n-k} : \, |z| < r}\subset \r
^n$, where $r>0$, and
$$D= \set{\left( \frac{s}{\sqrt{1+r^2+s^2}}
\theta _1 , \frac{t\sqrt{1+r^2}}{r\sqrt{1+r^2+s^2}} \theta _2
,\frac{\sqrt{1+r^2}\sqrt{r^2-t^2}}{r \sqrt{1+r^2+s^2} }\right) \in \s ^n :  \,
\begin{matrix}
\theta _1 &\in & \s ^{k-1}\\
\theta _2 &\in & \s ^{n-k-1}\\
s &\geq & 0 \\
t & < & r
\end{matrix}}.$$

Then,

\begin{eqnarray*}
\psi (x,z)&=& \left( \sqrt{|x|^2+1+r^2}, x , z , \sqrt{r^2 -|z|^2}\right)\\
\eta (x,z)&=& \left( \frac{r \sqrt{|x|^2 +1+r^2}}{\sqrt{1+r^2}}, \frac{r
x}{\sqrt{1+r^2}}, \frac{\sqrt{1+r^2}\, z}{r},
\frac{\sqrt{1+r^2}\sqrt{r^2-|z|^2}}{r}\right)
\end{eqnarray*}

Thus, from \eqref{psiG}, we get
\begin{eqnarray}
\rho (x,z) &=& \ln \left( \frac{(r+\sqrt{1+r^2})\sqrt{|x|^2 +1 +
r^2}}{\sqrt{1+r^2}}\right)\\
G(x,z) &=& \left( \frac{x}{\sqrt{|x|^2 +1 + r^2}}, \frac{\sqrt{1+r^2}\, z}{r
\sqrt{|x|^2 +1 + r^2}}, \frac{\sqrt{1+r^2}\sqrt{r^2 - |z|^2}}{r \sqrt{|x|^2 +1 +
r^2}}\right)
\end{eqnarray}

In this case, the hypersuface has two principal curvatures given by, $k_i
=-\frac{r}{\sqrt{1+r^2}} $, for $i =1, \ldots, k$, and $k_j =
-\frac{\sqrt{1+r^2}}{r}$, for $j=k+1, \ldots , n$. Thus, the eigenvalues of the
Shouten tensor associated to $g$ (given by \eqref{isop}) are $\lambda _i =
-\frac{1}{2} - r^2 + r \sqrt{1+r^2}$, for $i=1,\ldots ,k$, and $\lambda _ j =
\frac{1}{2} + r^2 - r \sqrt{1+r^2}$, for $j = k+1 , \ldots , n$.

\end{enumerate}

\begin{remark}
As we pointed out at the begin of the Section, hyperspheres and horospheres are the
only solutions for $\sigma _k (\lambda _i)= 1$ on $\s^n$ and $\sigma _k (\lambda _i)
=0$ on $\r^n$. The other cases define complete metrics on a \emph{subdomain} of the
sphere. So, the natural question is: \emph{Are these solutions the only solutions
for such domains under the constrain $\sigma _k (\lambda _i) = {\rm constant}$?}
\end{remark}

\section{Invariant conformal metrics on the sphere}

In this Section we will give sufficient conditions for a conformal metric on the
sphere to be radial. The following local result is based on the correspondence given
in Theorem \ref{representacion}, Lemma \ref{confinv} and a deep result of Do
Carmo-Dajzcer for hypersurfaces in hyperbolic space.

\begin{theorem}\label{t1}
Let $g= e^{2\rho} g_0$ be a conformal metric defined on a domain of the sphere
$\Omega \subset \s ^n$ such that the eigenvalues, $\lambda_i$, for $i =1, \ldots,
n$, of its Schouten tensor, ${\rm Sch}_g$, verify
$${\rm sup}\set{\lambda _i (x) , \, x \in \Omega , \, i = 1 ,\ldots , n }< +\infty .$$

Furthermore, assume that the eigenvalues satisfy
\begin{equation*}
\begin{split}
\lambda &= \lambda _1 = \cdots = \lambda _{n-1} \\
\nu = \nu (\lambda) &= \lambda _n \\
\lambda - \nu & \neq 0 .
\end{split}
\end{equation*}

Then, $g$ is radial.
\end{theorem}
\begin{proof}
Consider $t>0$ big enough such that the eigenvalues of the Schouten tensor of $g_t =
e^{2t}g$ are less than $1/2$ (see Remark \ref{remark}). Consider the horospherically
convex hypersurface, $\phi : \Omega \To \h ^{n+1}$, associated to $g_t$ given by
\eqref{repfor} in Theorem \ref{representacion}. Hence, the principal curvatures of
$\phi$ verify:
\begin{equation*}
\begin{split}
\tilde{\lambda} &= \kappa _1 = \cdots = \kappa _{n-1}\\
\tilde{\nu}=\tilde{\nu}(\tilde{\lambda}) &= \kappa _n \\
\tilde{\lambda} - \tilde{\nu} & \neq 0 ,
\end{split}
\end{equation*}this follows from \eqref{lambdakappa} and the assumptions on the
eigenvalues of ${\rm Sch}_g$. Hence, using \cite[Theorem 4.2]{CD}, $\phi (\Omega)$
is contained in a rotational hypersurface, which means, via Lemma \ref{confisoinv},
that $g_t$ is radial, so $g$ is radial.
\end{proof}

The next result is about determining which conformal metrics on the sphere are
invariant by a $k-$parameter subgroup of conformal diffeomorphisms of the sphere. We
should remove the class of conformal metrics on the sphere such that all the
eigenvalues of its Schouten tensor are constant, $\mathcal{C}(n)$, but this is not a
significant problem, since there are not too many of them and we have classify them.
Again, the result is based on a Theorem of M. Do Carmo and M. Dajczer.

\begin{theorem}\label{t2}
Let $g= e^{2\rho} g_0$ be a conformal metric defined on a domain of the sphere
$\Omega \subset \s ^n$ such that $g \not\in \mathcal{C}(n)$ and the eigenvalues,
$\lambda_i$, for $i =1, \ldots, n$, of its Schouten tensor, ${\rm Sch}_g$, verify
$${\rm sup}\set{\lambda _i (x) , \, x \in \Omega , \, i = 1 ,\ldots , n }< +\infty .$$

Suppose that $g$ is invariant by a $k-$parameter subgroup of conformal
diffeormorphism $\mathcal{G} \leq \Diff$. Then the maximum value of $k$ is $k_{max}=
\frac{n(n-1)}{2}$, and if $k = k_{max}$, the Schouten tensor of $g$, ${\rm Sch}_g$
has two eigenvalues $\lambda$ and $\nu$, where one of them, say $\lambda$, has
multiplicity at least $n-1$. If, in addition, $\lambda \neq 0$, $\nu = \nu
(\lambda)$ and $\nu - \lambda \neq 0$, then $g$ is radial.
\end{theorem}
\begin{proof}
As above, dilate $g$ until the eigenvalues of the Schouten tensor are less than
$1/2$. Now, construct the horospherically convex hypersurface given by Theorem
\ref{representacion}. The hypothesis on the $\mathcal{G}-$invariance of $g$ is
translated into a $\mathcal{T}-$invariance of $\phi$ under a $k-$parameter subgroup
$\mathcal{T}\leq \IsoH$. Thus, applying now \cite[Theorem 4.7]{CD} we obtain the
result.
\end{proof}

\begin{remark}
The above results hold for $n \geq 3$. It is clear that for $n=2$ are false.
\end{remark}

\section{A note on rotational hypersurfaces in $\h ^{n+1}$}

In a recent paper \cite{CHY}, authors have classified all possible radial solution
to the equation
$$ \sigma _k (\lambda _i) = c , \, \, c \equiv {\rm constant}, $$that is, they
consider conformal metrics $g = v(|x|)^{-2}|dx|^2$ on domains of the form
$$\set{x \in \r ^n , \, r_1 < |x| < r_2} ,$$ being $\sigma _k (\lambda _i)$ the $k-$th
elementary symmetric function of the eigenvalues of ${\rm Sch}_g$, and $0 \leq r_1 <
r_2 \leq \infty$.

From the point of view of hypersurfaces in hyperbolic space, this classification
result means (up to possibly a dilatation) that they have classified all rotational
horospherically convex hypersurfaces verifying the Weingarten relationship
$$ \sigma _k \left( \frac{1+\kappa _i}{2(1- \kappa _i)}\right) = \tilde{c} ,
\, \, \tilde{c} \equiv {\rm constant} .$$

It will be too long to describe here all these solutions, but we would like to
mention two cases when $c >0$: {\bf Case I.1} and {\bf Case I.3.a} in \cite[Theorem
1]{CHY} give the existence of hyperspheres (which was already known) and
Delaunay-type hypersurfaces respectively.

\begin{remark}
An interesting application of the above hypersurfaces could be to use them as
barriers for Plateau problem at infinity in the hyperbolic space for certain
Weingarten functionals.
\end{remark}

\vspace{.5cm}
\begin{center}
{\bf Acknowledgement:}
\end{center}
Author wants to thank J.A. Gálvez, Y.Y. Li and R. Mazzeo for their interesting
comments and help during the preparation of this work.

\end{document}